
\documentstyle{amsppt}
\magnification=\magstep1        
\NoRunningHeads
\TagsOnRight

\def\pint{\Cal{P}-\int}
\def\X{\frak X}
\def\lor{L_1 (\bold R)}
\def\lox{L_1 (\X)}
\def\pox{\Cal P_1 (\X)}
\def\rn{\right\Vert}
\def\ln{\left\Vert}
\def\lav{\left\vert}
\def\rav{\right\vert}
\def\e{\epsilon}
\def\qed{\hfill{\vrule height6pt  width6pt depth0pt}}
\font\cmcsc=cmcsc8
\def\D{\Delta}
\def\hf{\hfill}

\topmatter

\title
Nowhere Weak Differentiability of the Pettis Integral  
\endtitle
\date
Final Version:  31 October    1994
\enddate
\dedicatory
to appear:  Quaestiones Mathematicae
\enddedicatory
\author 
S.J. Dilworth   
\quad  {\cmcsc and}   \quad
Maria Girardi
\endauthor
\address
Department of Mathematics,
University of South Carolina,
Columbia, SC  29208, U.S.A.
\newline
Current address:
Department of Mathematics, Bowling Green State University,
Bowling Green, Ohio 43403, U.S.A.
\endaddress
\email dilworth\@math.scarolina.edu \endemail
\address
Department of Mathematics,
University of South Carolina,
Columbia, SC  29208, U.S.A.
\endaddress
\email
girardi\@math.scarolina.edu \endemail
\thanks
The second author was supported in part by NSF  DMS-9306460.
\endthanks
\abstract 
For an arbitrary infinite-dimensional Banach space $\X$,
we construct examples of  strongly-measurable   
$\X$-valued  Pettis integrable functions
whose indefinite Pettis integrals are nowhere weakly differentiable; 
thus, for these functions the 
Lebesgue Differentiation Theorem fails rather spectacularly.
We also relate the  degree of nondifferentiability 
of the indefinite Pettis integral to the cotype of $\X$, from which
it follows that  
our examples are reasonably sharp.  
\endabstract
\endtopmatter

\document
\baselineskip  14 pt

There are several generalizations of the
space $\lor$ of   Lebesgue integrable  functions
taking values in the  real numbers $\bold R$
(and defined on  the usual   
Lebesgue measure space $(\Omega, \Sigma, \mu)$ on $[0,1]$~)
to a
space of strongly-measurable ``integrable''  (suitably formulated)
functions  taking values in a Banach space  $\X$.

The most common  generalization is
the space $\lox$ of
Bochner-Lebesgue integrable functions.
Using the fact [P1, Theorem~1.1] that
a strongly-measurable    
function is
essentially separably-valued,
one can easily extend Lebesgue's Differentiation Theorem from
$\lor$ to  $\lox$.
Specifically  [B; cf\. DU, Theorem~II.2.9],
if $f\in \lox$, then
$$
\lim_{h\to 0} \
   \frac{1}{h} \int_t^{t+h}
   \ln f(\omega) - f(t) \rn \, d\mu(\omega)   ~=~  0
$$
and so
$$
\lim_{h\to 0} \
   \frac{1}{h} \int_t^{t+h} f(\omega) \, d\mu(\omega) ~=~  f(t)
$$
for almost all $t$ in  $\Omega$.

Another generalization of $\lor$  is the space $\pox$  of
strongly-measurable
Pettis integrable functions.
A  function $f : \Omega\to\X$
is Pettis integrable  if
for each $E\in\Sigma$ there is an element $x_E \in\X$
satisfying
$$
          x^*(x_E) ~=~ \int_E x^*f d\mu
$$
for each $x^*$ in the dual space $\X^*$ of $\X$.
The element $x_E$ is called the
Pettis integral of $f$ over $E$
and we write
$$
\pint_E f\, d\mu ~=~  x_E  \  .
$$ It is clear that
$\lox \subset \pox$, while the reverse inclusion holds if and only if
$\X$ is finite dimensional (see e\.g\.~[DG]).

If $f \in \pox$,
then  for each $x^*\in\X^*$
the function $x^* f \in \lor$
and so there exists a set $A$  (which depends on $x^*$)
of full measure  such that
$$
\lim_{h\to 0} \
   \frac{1}{h} \int_t^{t+h}  x^*f(\omega)\, d\mu(\omega)   ~=~  x^*f(t)
$$
for each  $t \in A$.
In  his paper  [P1] introducing the Pettis integral,
Pettis phrased this by saying that
the Pettis integral of a function in $\pox$ is
{\it pseudo-differentiable\,}.
He closed his paper by asking whether
the  Pettis integral of   a function $f$ in $\pox$
enjoys  the stronger property of being
{\it a\.e\.~weakly differentiable\,}; that is,  
does there necessarily exist a set $A$
({\it independent} of $x^*$) of full measure such that
$$
   \lim_{h\to 0} \
   \frac{1}{h} \int_t^{t+h}  x^*f(\omega)\, d\mu(\omega) ~=~  x^*f(t)
$$
for each  $t \in A$ and $x^* \in \X^*$, or such that (which is
the same thing of course)
$$
 \text{weak}-\lim_{h\to 0} \
    \frac{1}{h} \  \pint_t^{t+h} f(\omega)\, d\mu(\omega)  ~=~  f(t)
$$
for each  $t \in A$.

If $\X$ is finite dimensional,
then
the Pettis integral of a function  in $\pox$  is  
a\.e\.~weakly differentiable.   
R\.S\.~Phillips [Ph] (for $\X = \ell_2$)
and M\.E\.~Munroe  [M] (for $\X = C[0,1]$)
each constructed  an example of a
function  in $\pox$ whose Pettis integral is  not 
a\.e\.~weakly differentiable.  
G\.E\.F\.~Thomas [T, p\.~131] conjectured that 
such a function in $\pox$ exists
for every infinite-dimensional  Banach
space $\X$.

At the recent May 1993 Kent State University Functional Analysis
Conference,
Joe Diestel requested a further investigation 
into  Pettis's    
question.  
Independently, V\.~Kadets [K] recently constructed, 
for each infinite-dimensional Banach space $\X$,   
a function in $\pox$ whose Pettis integral fails    
to be a\.e\.~weakly differentiable; 
specifically,   it fails to be  
weakly differentiable on a  
set of positive, but not full, measure. 

The main theorem
of this paper  
constructs, 
for each infinite-dimensional Banach space $\X$,     
a function in $\pox$ whose Pettis integral is {\it nowhere}  
weakly differentiable.   
This theorem also  addresses the 
degree of nondifferentiability of the Pettis integral.   
Our second theorem shows, for {\it arbitrary} Banach spaces, that    
the functions which we construct are close to being 
optimal with respect to their  
degree of nondifferentiability.  From these  
two theorems it follows  
(Corollaries~3 and 4) that the     
cotype of a space is  closely tied to the degree of nondifferentiability
of the Pettis integral.

Theorem~2   was shown to us by Nigel
Kalton in answer to a question posed in a preliminary version of this paper.
We  are grateful   to him for permission to include this 
result here.

To state our main result we introduce
the  collection $\Psi$   of all  increasing   functions
$\psi \: \left[ 0, \infty \right) \to \left[ 0, \infty \right) $ satisfying
the growth condition
$$
 \sum_{n=1}^\infty
  \psi(2^{-p_{n-1}})~\sqrt{2^{p_n}}  ~<~ \infty\ ,
\tag"$(\dagger)$"
$$
for some increasing sequence $\{p_n\}_{n=0}^\infty$ of integers.
Examples of functions in $\Psi$ are
$$
\align
  &\quad\psi(s) = s^{^{\frac{1}{2} + \e}}, \\
\psi(s) = s^{\frac{1}{2}}\left[\frac{1}{\log\left(1/s\right)}\right]^{1+\e}
&\text{\qquad and \qquad}
\psi(s) = s^{\frac{1}{2}}  \left[ \frac{1}{\log\left(1/s\right)}\right]~
 \left[\frac{1}{\log\log\left(1/s\right)}\right]^{1+\e}
\endalign
$$
for $p_n = n$  and any $\e>0$.

\proclaim{Theorem~1}
Let $\X$ be an infinite-dimensional Banach space. For each
$\psi\in\Psi$,
there exists $f\in\pox$  such that
$$
 \ln \pint_I f\, d\mu \rn_{\X}  ~\geqslant~
\psi\left(\mu\left(I\right)\right) 
\tag"$(\ddagger)$"
$$
for each interval $I$ contained in \,$\left[0,1\right]$.  
\endproclaim
\remark{Remark}    
Taking $\psi(t) = t^{\frac{3}{4}}$    gives  a
Pettis integrable function $f$ such that for each $t \in \Omega$,
$$
\lim_{h\to 0}
  \ln \frac{1}{h} \  \pint_t^{t+h} f(\omega)\, d\mu(\omega) \rn_{\X}
   ~=~ \infty \ . 
$$
If the Pettis integral of this  $f$ were     
weakly differentiable at $t$,
then the above limit would be finite.
\endremark

\demo{Proof}
Let $\{ I^n_k : n=0,1, \ldots \ ; k=1, \ldots, 2^n \}$
be the dyadic intervals
on $[0,1]$, i\.e\.   
$$
  I^n_k  = \left[ \frac{k-1}{2^n}, \frac{k}{2^n} \right) \ .
$$
Define inductively a collection
$\{ A^n_k : n=0,1, \ldots \ ; k=1, \ldots, 2^n \} $ of disjoint sets
of strictly positive measure such that
$A^n_k \subset I^n_k$ (e\.g\.~appropriately chosen ``fat Cantor'' sets).

Fix $K > 1$. By a theorem of Mazur there is 
a basic sequence $\{ x_n \}$ in $\X$ with basis constant at
most $K$.
Take a blocking $\{ F_n \}$  of the  basis  with
each subspace $F_n$  of large enough dimension
to find (using the finite-dimensional version of
Dvoretzky's Theorem [D]) a $2^n$-dimensional subspace $E_n$ of $F_n$
such that the Banach-Mazur distance between
$E_n$ and $\ell_2^{2^n}$ is less than  2.
Note that  $\{ E_n \}$ forms a  finite-dimensional decomposition.
Next find   operators $T_n \: \ell_2^{2^n} \to E_n$
such that $\ln T_n \rn \leqslant 2$ and   $\ln T_n^{-1} \rn = 1$.
Let  $\{ u^n_k \: k=1,\ldots 2^n \}$ be the
standard unit vectors  of $\ell_2^{2^n}$
and let $  e^n_k  \equiv  T_n u^n_k$.

By the growth condition $(\dagger)$ on $\psi$, 
there is an  increasing sequence $\{p_n\}_{n=0}^\infty$ of integers,
with $p_0 = 0$, satisfying 
$$
\sum_{n=1}^\infty
  \psi(4~\cdot~2^{-p_{n-1}})~\sqrt{2^{p_n}}  ~<~ \infty \ . 
$$
Define $f \: \left[ 0,1 \right] \to \X$  by
$$
     f (\omega) ~=~
     \sum_{n=1}^\infty \sum_{k=1}^{2^n}
       ~c_n ~
      \frac{1_{A^n_k}(\omega)}{\mu(A^n_k)} ~  e^n_k  \  ,
$$
where
$$
   c_m ~=~  2K~
 \left[ \psi\left( 4~\cdot~2^{-{p_{n-1}}} \right) \right]
   \cdot~ \delta_{m,p_n} \  ,$$ 
(here $\delta_{j,k}$ is the usual Kronecker delta symbol).    
Clearly, $f$ is strongly measurable.

The Pettis integral of $f$ is easily computable; namely,
$$
\pint_E f\, d\mu ~=~
\sum_{n=1}^\infty \sum_{k=1}^{2^n}  ~c_n~
 \left(\int_E \frac{1_{A^n_k}}{\mu(A^n_k)} \, d\mu \right)
  ~ e^n_k \ .  
\tag"$(\ast)$"
$$ 
To see this, first note  that
the growth condition on $\psi$ guarantees that
the above series does indeed converge to an element
of $\X$, since
$$
\allowdisplaybreaks
\align
\ln  \sum_{n=p}^q  \sum_{k=1}^{2^n}  ~c_n ~
  \left(\int_E \frac{1_{A^n_k}}{\mu(A^n_k)} \, d\mu \right)
   e^n_k \rn_{\X}
~&=~
\ln  \sum_{n=p}^q  \sum_{k=1}^{2^n}  ~c_n ~
  \left(\int_E \frac{1_{A^n_k}}{\mu(A^n_k)} \, d\mu \right)
   T_n u^n_k \rn_{\X}
   \\
~&\leqslant~
2~   \sum_{n=p}^q ~c_n ~ \ln \sum_{k=1}^{2^n}
 \left(\int_E \frac{1_{A^n_k}}{\mu(A^n_k)} \, d\mu \right)  u^n_k
   \rn_{\ell_2^{2^n}}   \\
~&=~
2~\sum_{n=p}^q  ~c_n~  \left[
  \sum_{k=1}^{2^n}
  \left| \int_E  \frac{1_{A^n_k}}{\mu(A^n_k)} \, d\mu \right|^2
  \right]^\frac{1}{2} \\
&~\leqslant~
2~\sum_{n=p}^q c_n ~\sqrt{2^n} \ ,   
\endalign $$
which approaches zero as $p,q \rightarrow \infty$.
Now fix $E\in\Sigma$ and $x^*\in\X^*$   and
let $\e^n_k = \text{sign\,}(x^* e^n_k)$. Then
$$
\sum_{k=1}^{2^n}  \lav x^* e^n_k \rav
~=~
\lav \sum_{k=1}^{2^n}  \e^n_k x^*  T_n u^n_k  \rav
~\leqslant~
\ln T^*_n \rn ~ \ln x^* \rn ~
 \ln \sum_{k=1}^{2^n} \e^n_k u^n_k \rn_{\ell_2^{2^n}}
~\leqslant~ 2~ \ln x^* \rn ~ \left(\sqrt{2^n}\right) \ ,
$$
and so   
$$
\allowdisplaybreaks
\align
\int_E \  \sum_{n=1}^\infty
    \lav  \sum_{k=1}^{2^n}  ~c_n~
     \frac{1_{A^n_k}}{\mu(A^n_k)}~x^*(e^n_k) \rav \, d\mu
~&=~
\sum_{n=1}^\infty   \int_E
       \lav \sum_{k=1}^{2^n} ~c_n~
      \frac{1_{A^n_k}}{\mu(A^n_k)} ~x^*(e^n_k)
      \rav \, d\mu \\
~&\leqslant~
\sum_{n=1}^\infty  \sum_{k=1}^{2^n} ~c_n~
\left( \int_E  \frac{1_{A^n_k}}{\mu(A^n_k)}
      \, d\mu \right)
    \lav x^* e^n_k \rav  \\
~&\leqslant~
\sum_{n=1}^\infty \sum_{k=1}^{2^n} ~c_n~ \lav x^* e^n_k \rav \\
&\leqslant~
2~ \ln x^* \rn ~
\sum_{n=1}^\infty \ c_n  \left(\sqrt{2^n}\right) < \infty  \ .
\endalign
$$
Thus we may interchange the integral and summation below to see that
$$
\allowdisplaybreaks
\align
\int_E x^*f\, d\mu
~&=~
\int_E  \  \sum_{n=1}^\infty \sum_{k=1}^{2^n}  ~c_n~
     \frac{1_{A^n_k}}{\mu(A^n_k)} ~ x^*(e^n_k)  \, d\mu \\
~&=~
\sum_{n=1}^\infty  \int_E
       \sum_{k=1}^{2^n}   ~c_n~
     \frac{1_{A^n_k}}{\mu(A^n_k)} ~ x^*(e^n_k)  \, d\mu
~=~
x^*\left( \sum_{n=1}^\infty\sum_{k=1}^{2^n}  ~c_n~
\left(\int_E  \frac{1_{A^n_k}}{\mu(A^n_k)} \, d\mu\right)
   e^n_k \right)   \ ,
\endalign
$$
as needed for $(\ast)$.

Fix an interval $I\in\Sigma$.
Find a dyadic interval $I^m_j \subset I$ such that
$ 4~\mu (I^m_j) \geqslant \mu(I)$ and then
find $n$ such that $p_{n-1} \leqslant m < p_n$.
Let $P$ be the natural projection
from $\sum\oplus E_j$ onto $E_{p_n}$.
Since   $ \ln P \rn \leqslant 2 K$,
$$
\allowdisplaybreaks
\align
2 K~ \ln    \pint_I f\, d\mu \rn_{\X}
~&\geqslant~
\ln P\left(\pint_I f\, d\mu \right)\rn_{\X} \\
~&=~
c_{p_n}
\ln \sum_{k=1}^{2^{p_n}} \left(\int_I \frac{1_{A^{p_n}_k}}{\mu(A^{p_n}_k)}
   \, d\mu \right) e^{p_n}_k\rn_{\X} \\
~&\geqslant~
c_{p_n}
\ln \sum_{k=1}^{2^{p_n}}
\left(\int_I  \frac{1_{A^{p_n}_k}}{\mu(A^{p_n}_k)}
  \, d\mu \right) u^{p_n}_k\rn_{\ell^{2^{p_n}}_2} \\
~&=~
c_{p_n}
\left[ \sum_{k=1}^{2^{p_n}} \lav \int_I 
\frac{1_{A^{p_n}_k}}{\mu(A^{p_n}_k)}
    \, d\mu \rav^2 \right]^\frac{1}{2}    \ ,
\endalign
$$
and  so since $A^{p_n}_k \subset I^{p_n}_k \subset I^m_j \subset I$ for
some $k$,
$$
2 K~ \ln    \pint_I f\, d\mu \rn_{\X}
~\geqslant~
c_{p_n}
~~=~~
2 K ~ \psi\left( 4~\cdot~2^{-{p_{n-1}}} \right) \ .
$$
But $\psi$ is increasing and
$4~\cdot~2^{-{p_{n-1}}}~\geqslant~4~\cdot~2^{-m}~\geqslant~\mu(I)$
and so
$$
 \ln    \pint_I f\, d\mu \rn_{\X}
~\geqslant~ \psi\left( \mu\left(I\right) \right) \ .
$$
Thus $f$ satisfies the statement of the theorem.   \qed
\enddemo

The functions in $\Psi$ can be viewed as indicators 
of the degree of nondifferentiability  
(i\.e\.~ the poor ``averaging behavior'') 
of the indefinite Pettis integral. 
For instance, taking 
$$
\psi(s) =
s^{\frac{1}{2}}\left[\frac{1}{\log\left(1/s\right)}\right]^{1+\e},
$$ 
we deduce from Theorem~1 that there exists $f\in\pox$  such that,
not only do we have
$$
\lim_{h\to 0}  ~
  \ln \frac{1}{h} \  \pint_t^{t+h} f(\omega)\, d\mu(\omega) \rn_{\X}
   ~=~ \infty \ ,
$$
but even worse,
$$
\lim_{h\to 0} ~ h^{\frac{1}{2}} ~\cdot~
   \left[ \log{\left(\frac{1}{h}\right)} \right]^{1+\e}~
  \ln \frac{1}{h} \  \pint_t^{t+h} f(\omega)\, d\mu(\omega) \rn_{\X}
   ~=~ \infty \ 
$$
for  all $t \in \Omega$.

The next theorem shows that 
Theorem~1 comes
close to describing the {\it worst} type of averaging behavior of the
Pettis integral that can occur in an {\it arbitrary} infinite-dimensional
Banach space.  In particular, 
it shows that,  for spaces on which the identity operator is 
$(2,1)$-summing (i\.e\., ~spaces with the  Orlicz property),      
Theorem~1 fails to hold for the function $\psi(s) = s^{\frac{1}{2}}$.  
Thus,  the growth condition $(\dagger)$ on $\psi\in \Psi$ 
can {\it not} be replaced by  $\psi(s) = O(s^\frac{1}{2})$ as $s\to 0$. 
We do not know, however, whether it can be replaced 
by $\psi(s) = o(s^\frac{1}{2})$ as $s\to 0$.

\proclaim{Theorem 2}  
If the   identity operator on an infinite-dimensional  
Banach space $\X$ is $(q,1)$-summing    
for some $ 2 \leqslant q < \infty$,    
then, for every  $f\in\pox$, 
$$  
\ln \pint_t^{t+h} f\, d\mu \rn_{\X} ~=~ 
 o\left( h^{\frac{1}{q}}\right)
$$  
as $h \to 0^+$   
for $\mu$-a\.e\. $t$.  
\endproclaim 
\flushpar  
The proof below, which uses a factorization theorem of 
Pisier [P], was pointed out to us by Nigel Kalton.  
 
\demo{Proof}   
Fix $f\in\pox$ for an infinite-dimensional 
Banach space $\X$.     
Consider the operator 
$K\: L_\infty \to \X$ given by 
$$ K(g) = \pint_{\Omega} g(\omega) f(\omega) \, d\mu(\omega) \ . 
$$ 
We need  to show that 
$$
   \ln K \left( 1_{[0, t+h]} \right) - 
        K \left( 1_{[0, t]} \right) \rn_{\X}   
  ~=~  o\left( h^{\frac{1}{q}}\right) 
$$ 
as $h \to 0^+$ for $\mu$-a\.e\. $t$. 
Fix $\e >0$.

Since  $K$ is compact  
and since the dual of $L_\infty$   
has the approximation 
property,  
there is [e\.g\.~DU,~Thm\.~VIII.3.6]  
a decomposition $K = K_1 + K_2$, with 
$K_i \in \Cal L \left(L_\infty, \X \right)$,   
such that  
$K_1$ has finite rank and $K_2$ has norm  at most $\e^2$.  
It is enough to show that there is  
some constant $A$, which  depends only 
on $\X$ and $q$,  such that   for each $i$,  
$$
 \limsup_{h\to 0^+}~h^{-\frac{1}{q}} ~\ln K_i\left( 1_{[0, t+h]} \right) - 
        K_i\left( 1_{[0, t]} \right) \rn_{\X} 
~\leqslant ~ A ~ \e  \ ,    
\tag"$(\lozenge)$"
$$  
on a set of $\mu$-measure at least $1 - \e^q$. 

Towards this, consider [see e\.g\.~R] the natural surjective isometry 
$\tau \:  L_\infty \to C(\D)  $   
for the appropriate extremally disconnected compact  Hausdorff  
space $\D$.  Recall that 
$\tau$ takes an indicator function of a Borel set in $[0,1]$  to  an 
indicator function   
of a  clopen set in $\D$, say $\tau\left( 1_A \right) = 1_{\widehat A }$ 
in such a way that if $A\subset B \subset \Omega$, 
then    $\widehat A\subset \widehat B \subset \D$ 
and $\widehat{ B\setminus A} ~=~ \widehat B \setminus \widehat A$.   
Let $\widehat K_i$ be the composite map: 
$$  \widehat K_i \quad\: \quad
  C(\D) \qquad \longrightarrow^{^{\hskip-1.25em \tau^{-1}}}\hskip1em \qquad
 L_\infty \qquad\longrightarrow^{^{\hskip-1.25em K_i}}\hskip1em\qquad \X \ . 
$$

First we deal with $K_1$.  
We assume, without loss of generality, that $K_1$ is of rank one.
So the mapping $\widehat K_1$  
is of the form 
$$  
\widehat K_1
 \left( \varphi \right) ~=~ 
  \left[  \int_\D \varphi \, d\lambda  \right] ~x  
$$ 
for some norm one element $x$ in $\X$ and a  
finite regular signed Borel measure $\lambda$ on $\D$.
Thus 
$$
\align 
\ln  K_1 \left( 1_{\left[ 0, t+h \right]} \right) 
    -  K_1 \left( 1_{\left[ 0, t \right]} \right)  \rn_{\X} 
~&=~ 
\ln \widehat K_1 \left( 1_{\widehat{\left[ 0, t+h \right]}} \right) 
  - \widehat K_1 \left( 1_{\widehat{\left[ 0, t \right]}} \right) \rn_{\X} \\
~&=~ \lav \lambda \left(\widehat{\left[ 0,t+h \right]}\right) 
   - \lambda \left(\widehat{\left[ 0, t \right]}\right) \rav \\
~&=~  
\lav \alpha (t+h)  - \alpha (t) \rav  \ ,
\endalign 
$$ 
where  $\alpha \: [0,1] \to \Bbb R$ is given by  
$\alpha (t)  = \lambda \left( \widehat{\left[ 0,t \right]}\right)$.  
Since $\widehat{\left[ 0,t \right]} \subset \widehat{\left[ 0,t+h \right]}$ 
for positive $h$,  the function 
$\alpha$ is of bounded variation and so  
is differentiable $\mu$-almost everywhere. 
Thus,    
$\ln  K_1 \left( 1_{\left[ 0, t+h \right]} \right) 
    -  K_1 \left( 1_{\left[ 0, t \right]} \right)  \rn_{\X} ~=~ O(h)$   
$\mu$-a\.e\.~ 
and  so $(\lozenge)$ holds for any $q > 1$.

Now we deal with $K_2$.   
Fix $2 \leqslant q < \infty$.   
If the   identity operator on $\X$ is  $(q,1)$-summing,  
then  [P, Cor.~2.7] there is a probability measure $\nu$  on the    
Borel sets  of $\D$ such that 
the operator $\widehat K_2$ admits  a factorization of the form 
$$\vbox{
\settabs\+\indent 
&\hskip 1 true in &\hskip 1 true in &\hskip .7 true in& \cr
\+&\hf$C(\D)\quad$
  &\hf$\longrightarrow^{^{\hskip-1.25em \widehat K_2 }}\hskip1em$\hf
  &\hf$\X$\hf&\cr 
\+\cr
\+&\hf$\searrow_{^{\hskip-1.3em J}}$& 
     &$\nearrow_{^{\hskip-.1em T}}$\hf&\cr  
\+\cr
\+& &\hf$L_{q,1}(\nu)$\hf&\cr
}$$ 
where  
$J$ is the natural inclusion map and $T$ is a   
bounded linear operator with operator norm at most  
$C\|\widehat K_2\|\le C \e^2$, where $C$ depends only  on $\X$ and q.   
Here,  $L_{q,1}(\nu)$ is the usual  
Lorentz space of all real-valued $\nu$-measurable functions $f$ 
on $\D$ for which the norm $\ln f \rn_{q,1}$ is 
finite, where 
$$   
\ln f \rn_{q,1} ~=~ \int_0^\infty  t^{\frac{1}{q} - 1} f^*(t) \, dt\, 
$$ 
and $f^*$  is the non-increasing rearrangement of $\lav f \rav$.    
As above   
$$
\align 
\ln  K_2 \left( 1_{\left[ 0, t+h \right]} \right) 
    -  K_2 \left( 1_{\left[ 0, t \right]} \right)  \rn_{\X} 
~&=~  
\ln  K_2 \left( 1_{\left(  t, t+h \right]} \right) 
      \rn_{\X}  \\ 
~&=~ 
\ln \widehat K_2 \left( 1_{\widehat{\left( t, t+h \right]}} \right) \rn_\X \\
~&\leqslant~ C\e^2 \ln J \left(1_{\widehat{\left( t, t+h \right]}} 
   \right)\rn_{L_{q,1}\left(\nu\right)}  \ . 
\endalign 
$$   
Since  the non-increasing rearrangement of 
$ J \left(1_{\widehat{\left( t, t+h \right]}}\right)$ 
is just the indicator function  of the set 
$\left[0, \nu\left(\widehat{\left(t, t+h \right]}\right)\right)$,  we have 
$$
\ln J \left(1_{\widehat{\left( t, t+h \right]}} 
   \right)\rn_{L_{q,1}\left(\nu\right)} 
~=~ 
q~ \left[ \nu\left(\widehat{\left(t, t+h \right]}\right)  
   \right]^{\frac{1}{q}} \ , 
$$ 
and so 
$$   
h^{-\frac{1}{q}}
\ln  K_2 \left( 1_{\left[ 0, t+h \right]} \right) 
    -  K_2 \left( 1_{\left[ 0, t \right]} \right)  \rn_{\X} 
~\leqslant~  
C q\e^2 \left[ 
\frac{ \lav \beta\left( t+h \right) - \beta\left(t \right) \rav}{h}
       \right]^{\frac{1}{q}} \ ,  
$$
where $\beta\: [0,1] \to \Bbb R$ is given by  
$\beta (t)  = \nu \left( \widehat{\left[ 0,t \right]}\right)$. 
The function $\beta$ is increasing and hence differentiable $\mu$-almost 
everywhere.  Thus 
$$
 \limsup_{h\to 0^+}~ h^{-\frac{1}{q}}~ \ln K_2\left( 1_{[0, t+h]} \right) - 
        K_2\left( 1_{[0, t]} \right) \rn_{\X} 
~\leqslant~ C~ q ~ \e^2 ~ \left[ \beta^\prime(t) \right]^{\frac{1}{q}} 
$$ 
for $\mu$-a\.e\. $t$. 
>From  $\int_0^1 \beta^\prime(t) \, dt 
 \leqslant \beta(1) - \beta(0) \leqslant 1$, 
it follows that   
$\mu\left[ \beta^\prime(t)  \geqslant \e^{-q} \right] \leqslant \e^q$. 
Thus, on a set of measure at least $1 - \e^q$, 
$$
 \limsup_{h\to 0^+}~h^{-\frac{1}{q}} ~\ln K_2\left( 1_{[0, t+h]} \right) - 
        K_2\left( 1_{[0, t]} \right) \rn_{\X} 
~\leqslant ~ C ~q ~ \e  \ ,    
$$ 
which implies $(\lozenge)$ for $K_2$. 
\qed
\enddemo

Recall that the  identity operator on a space with finite cotype $q$ 
is $(q,1)$-summing.    Indeed, cotype plays a major r\^ole 
in the unfolding drama.   
To see this, consider  
a space $\X$  which 
contains a finite-dimensional decomposition 
$\dsize\sum\oplus E_n$  where  the Banach-Mazur distance between  
$E_n$ and $\ell_p^{2^n}$ is less than  $M$ for each $n$  for some fixed  
$1 \leqslant p \leqslant \infty$ and $M > 1$. 
By modifying Mazur's construction [see e\.~g\.~LT] of a basic   
sequence and using the fact  (a simple compactness 
argument suffices) that finite representability of 
$\ell_p$ is inherited by subspaces of 
finite codimension, 
it is possible to construct such   
a  finite-dimensional decomposition  in $\X$ 
whenever $\ell_p$ is finitely
representable in $\X$. 
By the Maurey-Pisier Theorem  \cite{MP},  
$\ell_{q_0}$ is  finitely representable in $\X$ 
where   
$2 \leqslant q_0 \leqslant \infty $ and 
$$
   q_0 ~=~ \inf \ \{ q \: \X \text{  has cotype } q\}  \ . 
$$ 
In the same spirit as in the proof of Theorem~1
(and with similar notation),  
for $1\leqslant p\leqslant \infty$    
let  $\Psi_p$ be the collection  
of all  increasing   functions
$\psi \: \left[ 0, \infty \right) \to \left[ 0, \infty \right) $ 
satisfying    
the growth condition  
$$
\sum_{n=1}^\infty
  \psi\left(2^{-p_{n-1}}\right)~ 
\left[2^{p_n}\right]^{\frac{1}{p}}~<~ \infty      
\tag"$(\dagger_p)$"
$$ 
for some increasing sequence $\{p_n\}_{n=0}^\infty$ of integers  
(following the convention that $1/\infty$ is $0$). 
For $1\leqslant p < \infty$,  a typical function in  $\Psi_p$ is  
$$
\psi(s) = s^{\frac{1}{p} + \e} 
$$  
with $p_n = n$ and  for any $\e>0$. 
For $p = \infty$,   $(\dagger_p)$ 
reduces  to 
the condition   
$$
 \lim_{s\to 0^+} \psi(s) = 0 \ . 
$$  
Fix $\psi\in \Psi_p$ and find 
an increasing sequence $\{p_n\}_{n=0}^\infty$ of integers,
with $p_0 = 0$, satisfying 
$$
\sum_{n=1}^\infty
  \psi(4~\cdot~2^{-p_{n-1}})~ \left[2^{p_n}\right]^{\frac{1}{p}}~<~ \infty \  
$$  
(again, $1/\infty$ is $0$). 
Define $f \: \left[ 0,1 \right] \to \X$  by
$$
     f (\omega) ~=~
     \sum_{n=1}^\infty \sum_{k=1}^{2^n}
       ~c_n ~
      \frac{1_{A^n_k}(\omega)}{\mu(A^n_k)} ~  e^n_k  \  ,
$$  
where  
$$
   c_m ~=~  2K~
 \left[ \psi\left( 4~\cdot~2^{-{p_{n-1}}} \right) \right]
   \cdot~ \delta_{m,p_n} \ , 
$$ 
where $K$ is the finite-dimensional decomposition  constant.  
Minor variations of the proof of Theorem~1
show that this  function $f$   satisfies  
$$
 \ln \pint_I f\, d\mu \rn_{\X}
  ~\geqslant~ \psi\left(\mu\left(I\right)\right)$$
for each interval $I$ contained in $[0,1]$.

Theorems 1 and 2, along with the above observations, 
give the following corollaries.  

\proclaim{Corollary 3}   
Let $\X$ be an infinite-dimensional Banach space with finite cotype and
let  
$q_0 ~=~ \inf\{ q \: \X_0 \text{  has cotype } q\}$. Then  
the following hold. 
\roster 
\item If $p >  q_0$, then  for each $f\in\pox$, we have 
      $$  
      \ln \pint_t^{t+h} f\, d\mu \rn_{\X} ~=~ 
      o\left( h^{\frac{1}{p}}\right)
      $$  
      as $h \to 0^+$ for $\mu$-a\.e\. $t$. 
\item If $p< q_0$, then there is an $f\in\pox$  such that  
      $$ \ln \pint_t^{t+h} f\, d\mu \rn_{\X}  ~\geqslant~ 
        h^{\frac{1}{p}} 
      $$  
      for all $t \in [0,1]$.   
\endroster
\endproclaim

\proclaim{Corollary 4}  
For an  infinite-dimensional Banach space~$\X$,     
the following are equivalent.  
\roster
\item $\X$  fails cotype. 
\item For each $\psi \in \Psi_\infty$,  there exists $f \in \pox$ such that
$$
 \ln \pint_I f\, d\mu \rn_{\X}
  ~\geqslant~ \psi\left(\mu\left(I\right)\right)$$
for each interval $I$ contained in $[0,1]$.
\endroster
\endproclaim

\remark{Remark} Note that Corollary~4 proves the existence
of  a {\it reflexive} Banach space  for which the Pettis integral
has essentially no kind of
differentiability  property whatsoever.
\endremark

Theorem~1 can be reformulated by considering the indefinite
Pettis integral
$$g(t)=\pint_0^t f(\omega)\, d\mu(\omega) \ ,$$
and  then  expressing $(\ddagger)$ as 
$$
\ln g(s)-g(t) \rn \geqslant \psi(\lav s-t \rav )  \ . 
\tag"$(\ddagger^\prime)$"
$$ 
Corollary~4  
shows that if $g$ is the indefinite integral of  a  
Pettis-integrable function taking values in 
a space failing cotype, then there are (essentially) no restrictions 
on $\psi$ in $(\ddagger^\prime)$.  
Since $g(t)$ is always {\it continuous\,} 
[P1, Thm\.~2.5],    
it is not unreasonable to inquire,  
in the case of an arbitrary infinite-dimensional
Banach space, whether there are any 
restrictions on $\psi$ which are  attributable merely to the
continuity of $g$ as opposed to the additional fact that $g$ is an
indefinite
Pettis integral. 
Our final result answers this
question with a resounding no.  

\proclaim{Theorem~5}
Let $\X$ be an infinite-dimensional Banach space 
and let $\psi\in\Psi_\infty$. Then
there exists
a continuous function
$ f \: \Omega \to \X$ such that
$$
\ln   f(s) -  f(t) \rn_{\X}  
~\geqslant~
\psi\left(\lav s-t \rav  \right)
$$
for each $s$ and $t$ in $\Omega$.
\endproclaim
\remark{Remark}  
As Ralph Howard pointed out,  
Theorem~5 does not hold if $\X$ is finite-dimensional. In fact,
if $f$ is a
continuous function taking values in $\Bbb R^n$ and satisfying the
lower estimate given above, then 
an easy Hausdorff dimension argument
(see e\.g\.~[Kah]) shows that  the function $\psi$
must satisfy $\liminf_{t\rightarrow 0} \psi(t)t^{\e -1/n} < \infty$
for every $\e>0$. 
\endremark

\demo{Proof}
Find an increasing sequence $\{p_n\}_{n=0}^\infty$ of integers
with $p_0 = 0$ such that $\sum_n \psi \left(2^{-p_n}\right)$ is finite
and fix $K >1$.
Keeping with the notations and ideas of Theorem~1,
find a
finite-dimensional decomposition $\{ E_n\}$
in $\X$
and, to avoid excessive superscripts,
let  $J^n_k = I^{p_n}_k$ and likewise $\tilde e^n_k = e^{p_n}_k$
and $\tilde u^n_k = u^{p_n}_k$
for each admissible $n$ and $k$.

Consider the continuous
piecewise-linear function
$$
  f_n \left(\omega\right) ~=~
   \sum_{k=1}^{2^{p_n}}  ~ 2^{p_n} ~
\left[ \left(\frac{k}{2^{p_n}} - \omega\right)  ~ \tilde e^n_k
 + \left( \omega - \frac{k-1}{2^{p_n}}\right) ~ \tilde e^n_{k+1}
 \right] 1_{J^n_k} \left(\omega\right)  .
$$
If  $\omega\in J^n_k$, then $f_n\left(\omega\right)$ is of
the form $\alpha ~ \tilde e^n_k ~+~ (1-\alpha) ~\tilde e^n_{k+1}$
for some $0 \leqslant \alpha \leqslant 1$.
Thus  the norm of
$f_n\left(\omega\right)$ is  at most 2
for each $\omega\in\Omega$.
Define $f\: \Omega \to \X$  by
$$
f \left(\omega\right) ~=~
     \sum_{n=2}^\infty  ~c_n ~ f_n\left(\omega\right) \   ,
$$
where
$$
   c_{n+2} ~=~  2~K~ \psi\left( 2^{-p_{n}} \right) \ .
$$
Since each $f_n$ is uniformly continuous and
$$
\ln \sum_{n=p}^q c_n  f_n \left(\omega\right) \rn
~\leqslant 2~\sum_{n=p}^q  c_n      \ ,
$$
the choice of $\{ p_n \}$ guarantees  not only that
$f(\omega)$ is indeed in $\X$
for each $\omega\in\Omega$
but also that $f$ is uniformly  continuous.

Fix $s,t\in\Omega$.
Find $p_n$ such that $2^{-p_n} < \lav s-t \rav \leqslant 2^{-p_{n-1}}$.
Since $s$ and $t$ are in  neither the same nor adjacent
intervals of the partition $\{ J^{n+1}_k \}_k$ of $\Omega$,
for appropriate {\it distinct}
integers $k-1$, $k$, $j$, and  $j+1$,
$$
\align
f_{n+1}\left(s\right) ~&=~ \
 \alpha ~ \tilde e^{n+1}_{k-1} ~+~ (1-\alpha) ~\tilde e^{n+1}_k \\
f_{n+1}\left(t\right) ~&=~ \
 \beta ~ \tilde e^{n+1}_{j} ~+~ (1-\beta)~ \tilde e^{n+1}_{j+1} \,
\endalign
$$
for some  $0 \leqslant \alpha, \beta \leqslant 1$  and so
$$
\allowdisplaybreaks 
\align 
\ln f_{n+1}\left(s\right) - f_{n+1}\left(t\right) \rn_{\X}
~&\geqslant~
\ln  \alpha~\tilde u^{n+1}_{k-1} + (1-\alpha)~\tilde u^{n+1}_k
~-~  \beta~\tilde u^{n+1}_{j} - (1-\beta)~\tilde u^{n+1}_{j+1} \rn_{\ell_2} \\
&=~ \left[ \left(\alpha\right)^2 ~+~ \left( 1 - \alpha\right)^2
 ~+~ \left(\beta\right)^2 ~+~ \left(1-\beta\right)^2 \right]^{\frac{1}{2}}
\\
&\geqslant \quad
1 \ .
\endalign
$$
Let $P$ be the natural projection
from $\sum\oplus E_j$ onto $E_{p_{n+1}}$.
Since $\psi$ is increasing, we see that
$$
\align
2~K~\ln f\left(s\right) - f\left(t\right) \rn_{\X}
~&\geqslant~
\ln  P \left( f\left(s\right)
            - f\left(t\right) \right) \rn_{\X}\\
~&=~
c_{n+1}~\ln \left( f_{n+1} \left(s\right) -
                   f_{n+1} \left(t\right)\right)\rn_{\X} \\
~&\geqslant~
c_{n+1} \\
~&=~
2~K~\psi\left(2^{-p_{n-1}}\right) \\
~&\geqslant~
2~K~\psi\left(\lav s-t \rav \right) \ .
\endalign
$$
Thus $f$ satisfies the statement of the theorem.   \qed
\enddemo
\remark{Remark} Theorem~5 really only uses the existence of a basic
sequence inside $\X$, while Theorem~1 makes full use of Dvoretzky's
Theorem. 
\endremark
We close with a few observations.
[DG, Ex\.~3] constructs,
for each fixed infinite-dimensional Banach space~$\X$,
a strongly-measurable 
$\X$-valued function that is
Pettis integrable but not Bochner-Lebesgue integrable; however,
that function {\it is\,} Bochner-Lebesgue integrable over {\it any\,}
interval not containing 0.
Theorem~1 pushes this construction a bit further to give a
Pettis integrable function that  {\it is not\,}   
Bochner-Lebesgue integrable over  
{\it  any\,} interval.  

Consider the collection $K(\mu,\X)$
of the $\mu$-continuous countably additive  $\X$-valued vector measure
with relatively compact range.
If $f$ is in $\pox$, then the corresponding measure
$ \nu_f (E) = \pint_E f \, d\mu $
is in  $K(\mu,\X)$
[cf\.~DU, Thm\.~VIII\.1\.5].
The measure $\nu_f (E)$ is of bounded semi-variation;
furthermore,  $\nu_f (E)$ is of bounded variation
if and only if $f$ is in $\lox$
[cf\.~DU, Thm\.~II\.2\.4, Cor\.~2\.5].
Theorem~1
(consider the measure $\nu_f$ corresponding to $f$ as above)
and [JK, Theorem~2]
both  construct,
for each fixed infinite-dimensional Banach space $\X$,
a vector measure in  $K(\mu,\X)$ that
is of   bounded semi-variation but of
infinite variation on every interval.
The measure in [JK, Theorem~2] cannot arise, however, as an
indefinite Pettis integral, while the measure  from Theorem~1
is (of course) precisely an indefinite Pettis integral.

\widestnumber\no{[GGG]Z}
\def\n #1{\no{[\bf #1]}}

\enddocument